\documentclass{amsart}    
\usepackage{amsmath, amsthm, amssymb}
\usepackage[english]{babel}
\usepackage[T1]{fontenc}
\usepackage[utf8]{inputenc}
\usepackage{mathrsfs}
\usepackage{pgf}
\usepackage{pinlabel}

\title{An infinite clique of high-filling rays\\
 in the plane minus a Cantor set}

\author{Juliette Bavard}
\address{Univ Rennes, CNRS, IRMAR - UMR 6625, F-35000 Rennes, France}
\email{juliette.bavard@univ-rennes1.fr}

\newtheorem*{theo*}{Theorem}
\newtheorem*{prop*}{Proposition}

\theoremstyle{definition}

\newcommand{\N}{\mathbb{N}}
\newcommand{\R}{\mathbb{R}}

\DeclareMathOperator{\Cantor}{Cantor}

\begin{document}

\maketitle

\begin{abstract} The study of the mapping class group of the plane minus a Cantor set uses a graph of loops, which is an analogous of the curve graph in the study of mapping class groups of compact surfaces. The Gromov boundary of this loop graph can be described in terms of "cliques of high-filling rays": high-filling rays are simple geodesics of the surface which are complicated enough to be infinitely far away from any loop in the graph. Moreover, these rays are arranged in cliques: any two high-filling rays which are both disjoint from a third one are necessarily mutually disjoint. Every such clique is a point of the Gromov boundary of the loop graph. Some examples of cliques with any finite number of high-filling rays are already known. 

In this paper, we construct an infinite clique of high-filling rays.  \end{abstract}

\section{Introduction}

Because it has infinite topological type, the plane minus a Cantor set is quite richer than finite type surfaces. In particular, its uncountable mapping class group, as well as those of other infinite type surfaces, known as \textit{big mapping class groups}, let us with a whole new range of new phenomena to be described and understood. 

One way to progress in the study of these big mapping class groups is to make them act by isometries and non trivially -- when it is possible -- on Gromov-hyperbolic graphs. Following an idea described by Danny Calegari in \cite{Calegari-blog}, we defined in \cite{Juliette} the loop graph. The vertices of the graph are the simple loops on the surface, where a loop is an oriented geodesic ray starting and ending at infinity. An edge is defined as a couple of two disjoint loops. Each edge has length one. We proved in \cite{Juliette} that this loop graph has infinite diameter and is Gromov-hyperbolic. In \cite{Bavard-Walker}, together with Alden Walker we described the Gromov boundary of the loop graph in terms of cliques of high-filling rays. More precisely, we defined an extension of the loop graph by adding as vertices all the simple oriented geodesic rays of the surface which start at infinity. Edges are defined as before between any two disjoint rays. It turns out that this extended graph is no longer connected: it has infinitely many connected components. One of its component contains the loop graph, and is in fact quasi-isometric to it. In particular, it has infinite diameter and is Gromov-hyperbolic. All other components are cliques: they have diameter $0$ or $1$. The vertices in these components are called \textit{high-filling rays}: they are infinitely far away from any loop.

These high-filling rays are very interesting for many reasons. First, note that being \textit{filling} in the usual sense is not enough to be high-filling. Indeed, Lvzhou Chen and Alexander Rasmussen recently constructed in \cite{Chen-Rasmussen} uncountably many rays on the surface which intersect all loops but are at distance $2$ from some loop. In other words, there is another ray on the surface which is disjoint from both such a filling ray and a loop. These filling rays are in the component of the loop graph in the extended loop graph, and thus are not high-filling.

 Second, the Gromov-boundary of the loop graph can be identified with the set of cliques of high-filling rays (see \cite{Bavard-Walker}). Moreover, the two attractive and repulsive points of the boundary which are fixed by a loxodromic element of the mapping class group are necessarily \textit{finite} cliques, and have the same cardinal (see \cite{Bavard-Walker2}). On finite type surfaces, these high-filling rays would roughly correspond to the singular leaves of a minimal foliation. Using this correspondence and taking finite type subsurfaces of the plane minus a Cantor set give us cliques with any finite number of high-filling rays (see \cite{Bavard-Walker2}). But because this construction uses only \textit{finite} type subsurfaces, these cliques are not specific to infinite type surfaces. However, in \cite{Morales-Valdez}, Israel Morales and Ferr\'an Valdez constructed loxodromic elements which do not preserve any finite type subsurface and which have any chosen finite number of high-filling rays in their cliques. These high-filling rays are thus specific to infinite type surfaces. 

As far as we know, the question of whether or not there exist cliques containing infinitely many high-filling rays remains open. Note that the rays in such a clique would necessarily be specific to infinite type surfaces.  In this paper, we construct such an infinite clique of high-filling rays in the plane minus a Cantor set: more precisely, we obtain a clique with uncountably many high-filling rays.

\subsection*{Acknowledgment} We thank Jeremy Kahn and Alden Walker for helpful discussions on high-filling rays. We acknowledge support from the Centre Henri Lebesgue ANR11-LABX-0020-01.

\section{Dyadic tree and notations}

\subsection{General idea}
Our goal will be to embed a dyadic rooted tree in $\R^2 - \Cantor$, with the root at infinity, in such a way so that each branch of the tree is a high-filling ray. More precisely, any ray which is disjoint from the embedded tree will be isotopic to one of the branches of the tree. We will construct such an embedding by induction.

\subsection{Levels of the dyadic tree}
We consider the following tree:

\begin{figure}[!h]
\labellist
\small\hair 2pt
\endlabellist
\centering
\includegraphics[scale=0.8]{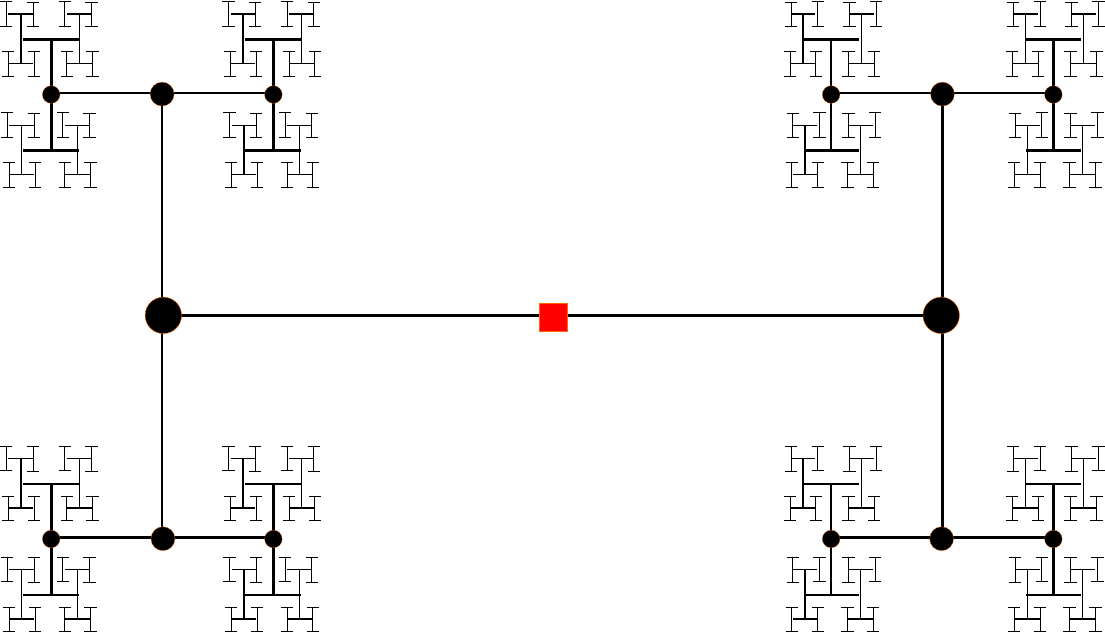}
\caption{Dyadic tree (the red squared dot is the root)}
\label{figu:tree1}
\end{figure}

We put a metric on this tree by deciding that each edge has length $1$. In particular, for all $n \in \N$, there are $2^n$ vertices at distance $n$ from the root. We call \emph{$n^{th}$ level of the tree} the subtree which contains all the vertices at distance at most $n$ from the root. We call \emph{branch} any geodesic path (finite or infinite) of the tree which starts at the root.

\begin{figure}[!h]
\labellist
\small\hair 2pt
\endlabellist
\centering
\includegraphics[scale=0.8]{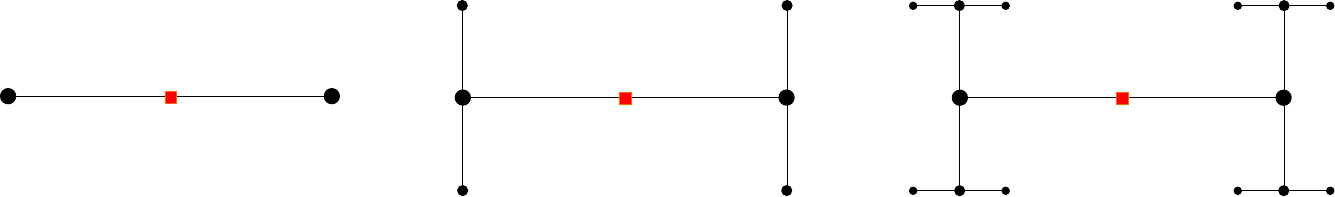}
\caption{First, second and third levels of the tree}
\label{figu:tree2}
\end{figure}

We will embed the tree by induction, level after level. We first give an embedding of the first and the second levels. We then explain how to extend the embedding of the $k^{th}$ level to the embedding of the $(k+1)^{th}$ level.

\subsection{Nested disks}

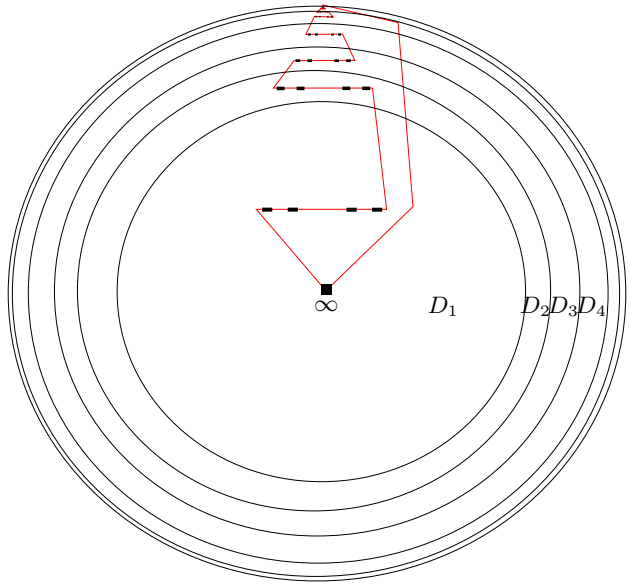
\begin{figure}[!h]
\labellist
\small\hair 2pt
\pinlabel $\infty$ at 273 236
\pinlabel $D_1$ at 373 236
\pinlabel $D_2$ at 452 236
\pinlabel $D_3$ at 476 236
\pinlabel $D_4$ at 500 236
\endlabellist
\centering
\includegraphics[scale=0.45]{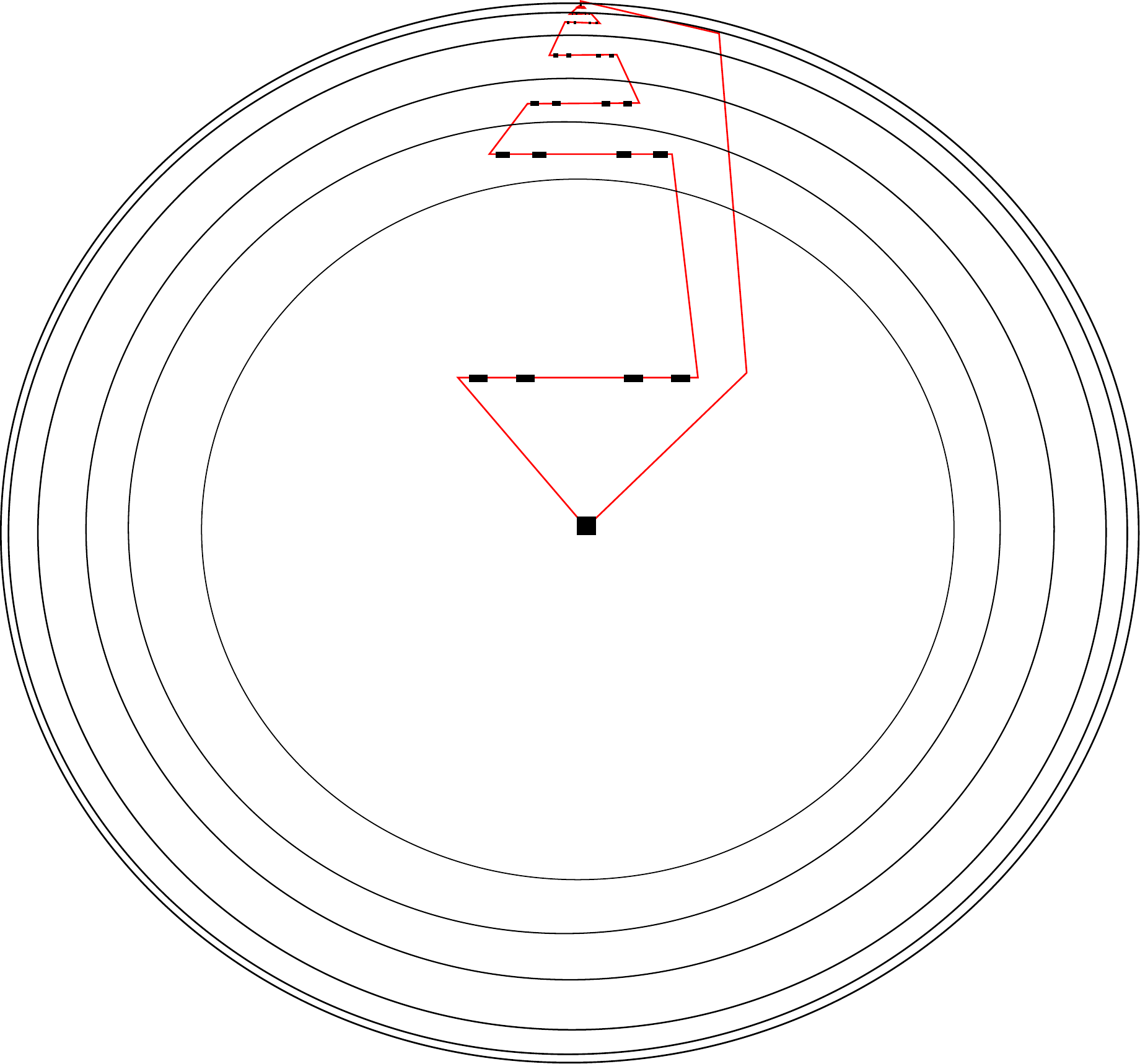}
\caption{Sequence of nested disks $(D_n)$ and an equator (in red)}
\label{figu:equator}
\end{figure}

We fix a complete hyperbolic metric of the first kind on $\R^2 - \Cantor$  (see for example the first part of \cite{Bavard-Walker2} to construct such a metric). We see $\R^2 - \Cantor$ as a sphere minus $\{\infty\} \cup \Cantor$ and speak about $\infty$ and the Cantor set as if they were marked points of the sphere. Choose a sequence of closed nested disks $(D_n)$  on the sphere such that $D_1$ contains $\{\infty\}$ and some points of the Cantor set, and such for all $i$, $D_{i+1}$ is not isotopic to $D_i$ in $\R^2 - \Cantor$ (i.e. $D_{i+1} - D_i$ contains points of the Cantor set). Moreover, we assume that only one point of the Cantor set is not included in the union of the $D_n$'s. This point is then in the closure of the union of the disks $D_n$'s. See Figure \ref{figu:equator} for an example of such a sequence of disk.

  We choose an equator of the sphere, i.e. a topological circle which contains $\{\infty\} \cup \Cantor$. See for example the red circle of Figure \ref{figu:equator}. For all $k \in \N$, the $k^{th}$ level of the tree will be embedded in $D_{k}$.

\subsection{Coding rays} Assume the boundaries of the $D_k$'s and the equator are geodesics. Following \cite{Juliette}, we can code the rays in the following way:
\begin{itemize}
\item We call \emph{equatorial segment} any subarc of the equator which has two endpoints in $\{\infty\} \cup \Cantor$ (in red on Figure \ref{figu:equator}). The equator is a countable union of its segments.
\item Any long ray is uniquely determined by the sequence of equator segments that its geodesic representative intersects from infinity (and the first hemisphere that it crosses -- north or south).
\end{itemize}

\section{Embeddings of the two first levels of the tree}
We will now embed the first and second levels of the dyadic tree in $D_1$ and $D_2$ respectively.

\subsection{First level}
\begin{figure}[!h]
\labellist
\small\hair 2pt
\pinlabel $\infty$ at 273 236
\pinlabel $D_1$ at 373 236
\endlabellist
\centering
\includegraphics[scale=0.45]{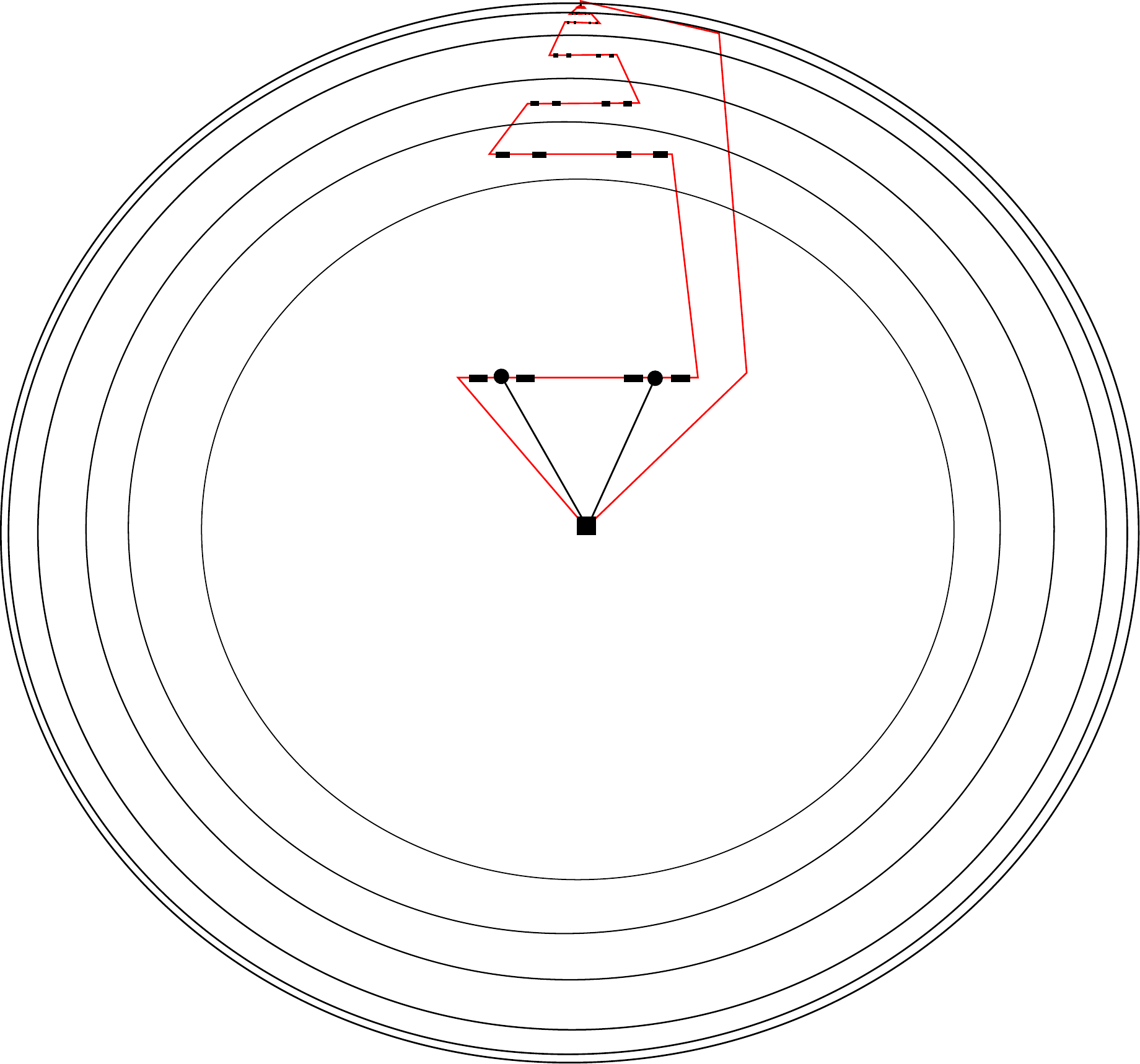}
\caption{Embedding of the first level}
\label{figu:level1}
\end{figure}

We embed the first level of the tree as in Figure \ref{figu:level1}: the root is at infinity, and the two branches (of length $1$) are included in the same hemisphere. They do not intersect the equator except at their ends. These ends lie in two different equatorial segments, in between points of the Cantor set included in $D_1$.

\subsection{Second level}
Here we explain how to extend the embedding of the first level of the tree to the second level. The method will be generalized in the next section.

We have embedded the first level in such a way so that:
\begin{itemize}
\item It is included in $D_1$.
\item The vertices at distance $1$ from the root are included in the equator.
\item We have two chunks of Cantor set on the left and of the right of each vertex on the equator.
\end{itemize}

Assume we have an embedding of the tree such that each vertex is included in an equator segment. We say that a ray \emph{begins like} some finite sub-branch of the tree if it crosses the same segments of the equator, in the same order, until it crosses the equatorial segment where the branch ends. 

We say that a ray or a branch \emph{begins like the $k^{th}$ level of the tree} if it begins like one of the $2^k$ branches of length $k$.

\begin{figure}[!h]
\labellist
\small\hair 2pt
\endlabellist
\centering
\includegraphics[scale=0.5]{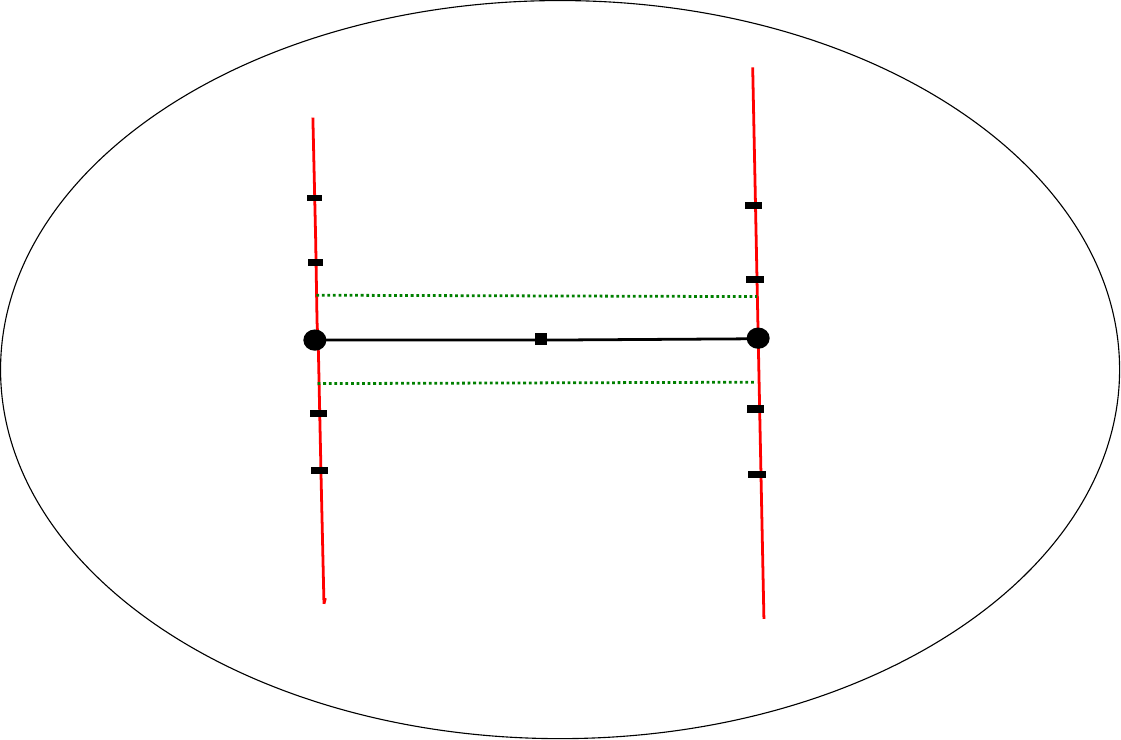}
\caption{Condition to begins like the first level}
\label{figu:level2-2}
\end{figure}

We want to embed the second level in such a way so that any ray disjoint from it begins like a branch of the first level. This will happen if and only if all the branches contain the two dotted green segments of Figure \ref{figu:level2-2}.

\begin{figure}[!h]
\labellist
\small\hair 2pt
\endlabellist
\centering
\includegraphics[scale=0.5]{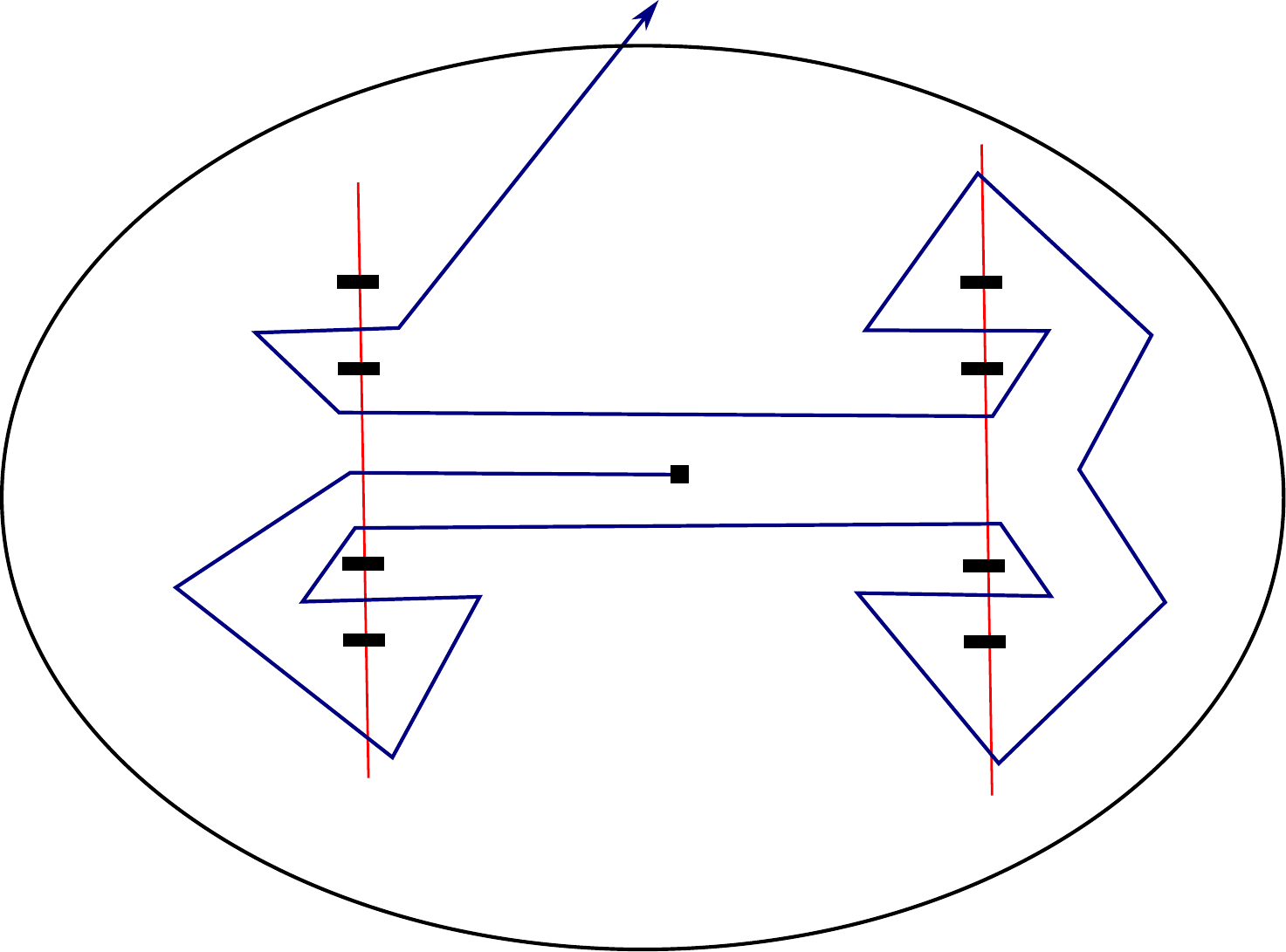}
\caption{First branch of level $2$ (only the part whichi is in $D_1$)}
\label{figu:level2-3}
\end{figure}

 We have $4$ branches to embed. We embed the beginning of a first one as in Figure \ref{figu:level2-3} (until it leaves $D_1$). Here we see that any ray disjoint from the blue branch of Figure \ref{figu:level2-3} will begin like the first level. Indeed, the drawing satisfies the condition of Figure \ref{figu:level2-2}.

\begin{figure}[!h]
\labellist
\small\hair 2pt
\endlabellist
\centering
\includegraphics[scale=0.5]{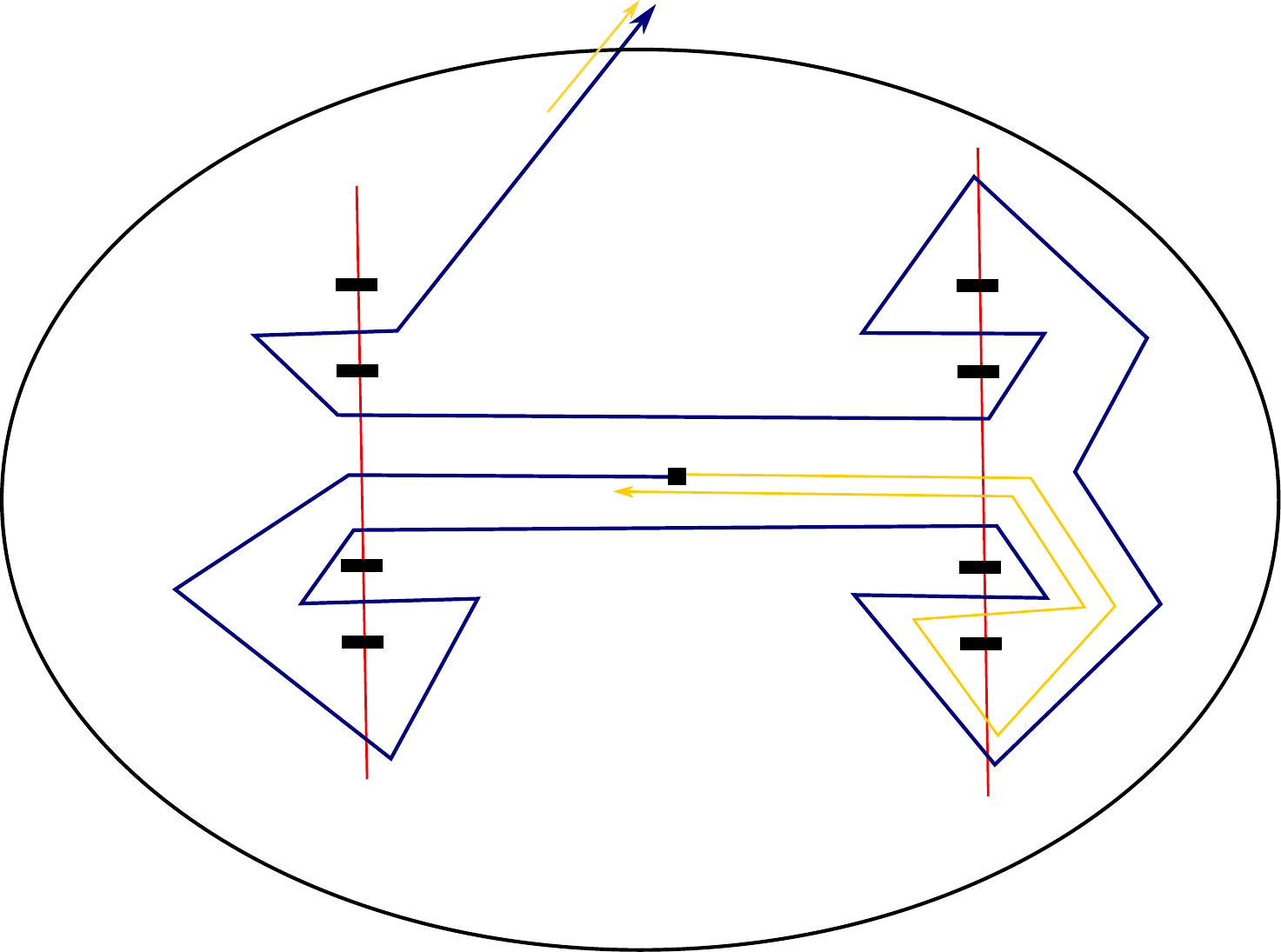}
\caption{Second branch of level $2$}
\label{figu:level2-4}
\end{figure}

We then embed the second branch as in Figure \ref{figu:level2-4}: it follows the yellow path and as soon as it comes back to infinity, it follows the first branch by staying in a tubular neighborhood on its left. In that way we are sure that the second branch will also satisfy the condition of Figure \ref{figu:level2-2} (i.e. it contains the two green dotted lines). Moreover, it leaves the disk $D_1$ on the tubular neighborhood of the first branch, on its left.

\begin{figure}[!h]
\labellist
\small\hair 2pt
\endlabellist
\centering
\includegraphics[scale=0.5]{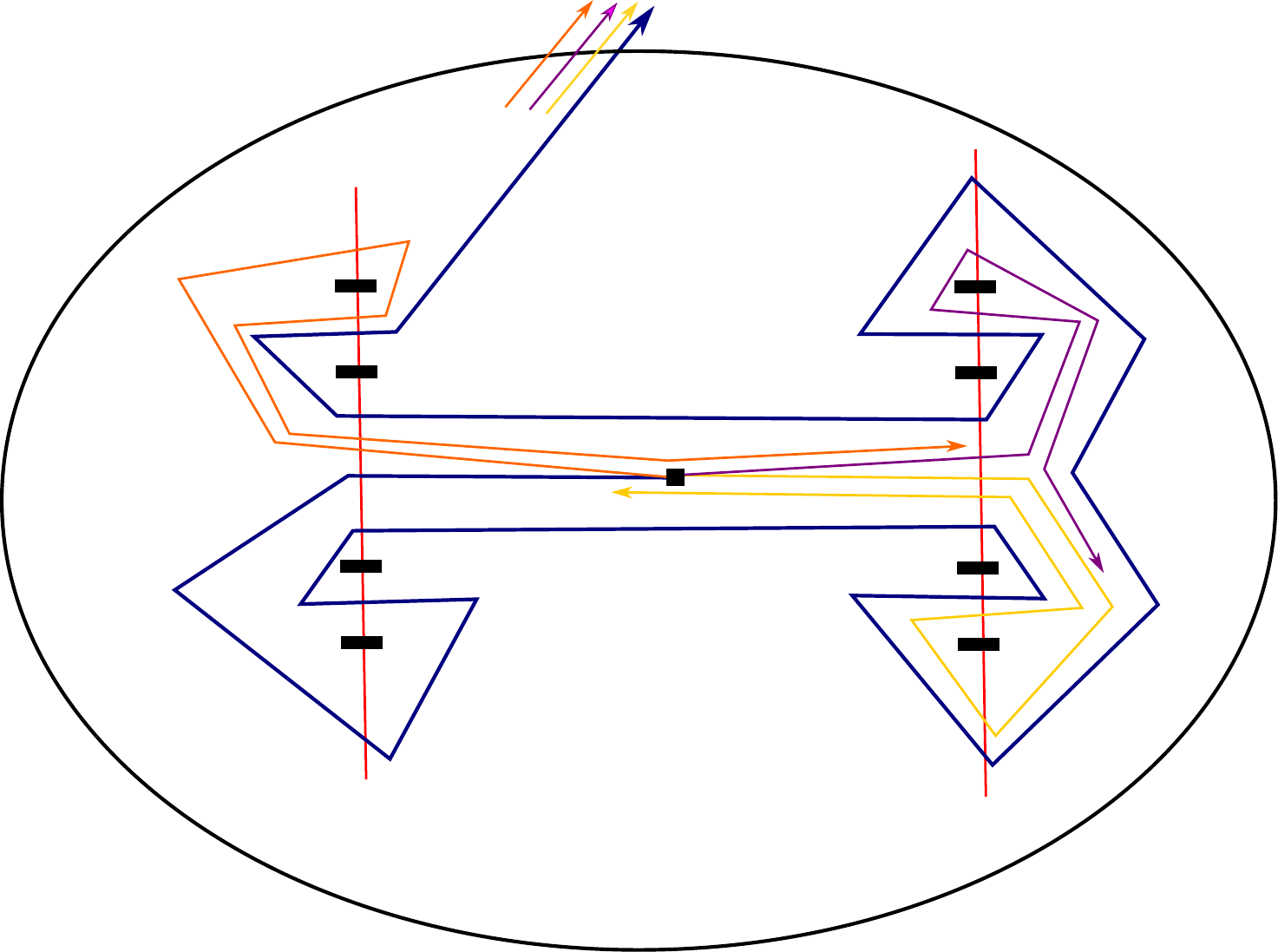}
\caption{Third and forth branches of level $2$}
\label{figu:level2-5}
\end{figure}

\begin{figure}[!h]
\labellist
\small\hair 2pt
\pinlabel $D_1$ at 343 136
\pinlabel $D_2$ at 376 136
\endlabellist
\centering
\includegraphics[scale=0.8]{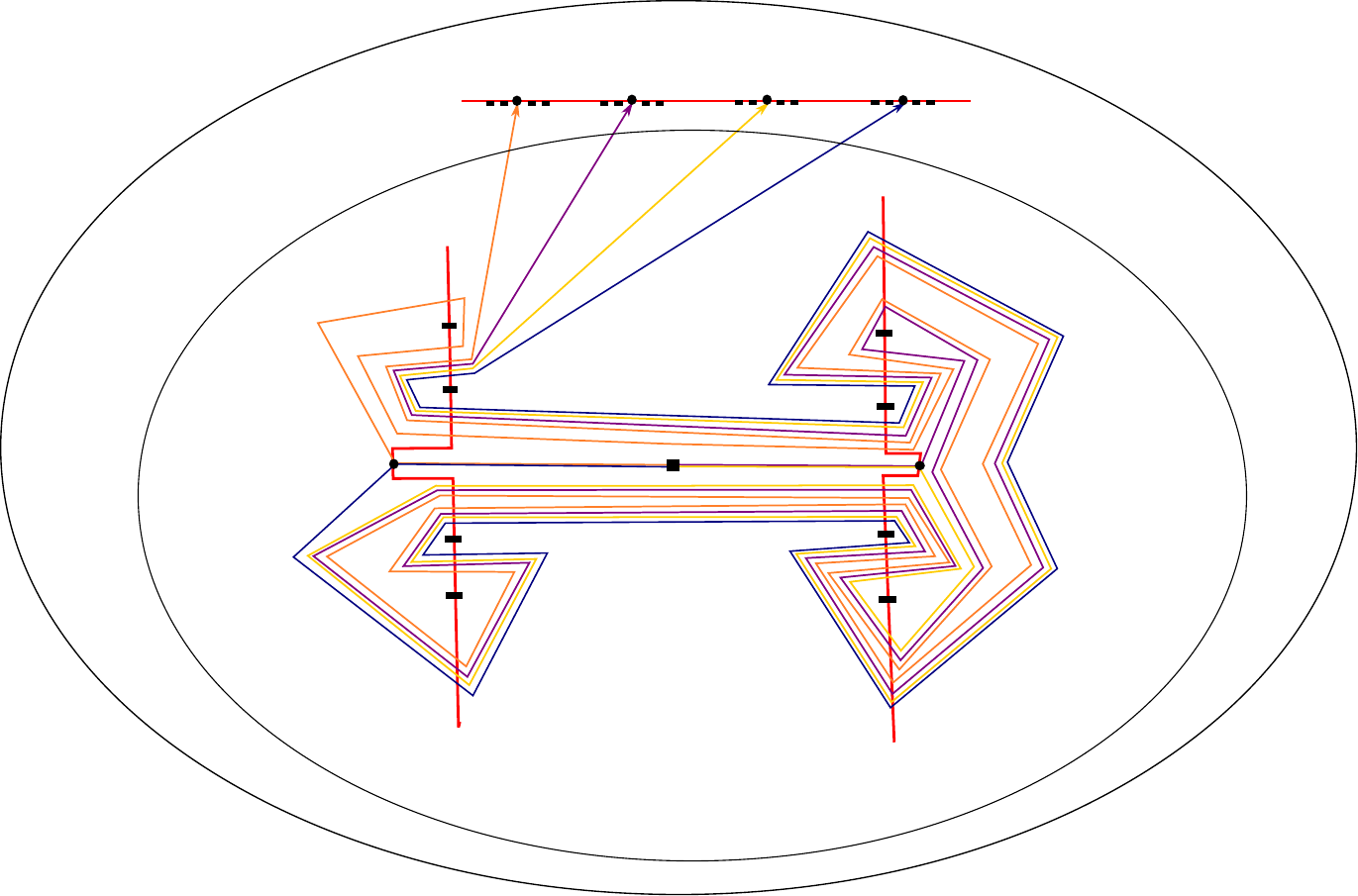}
\caption{Embedding of the second level (equator is in red)}
\label{figu:level2}
\end{figure}

We use the same idea to embed the third, and then the forth branches as in Figure \ref{figu:level2-5}: the third one follows the purple path and as soon as it meets the second branch (yellow), it follows it by staying in a tubular neighborhood on its left. In that way we are sure that the third branch will also follow the first branch, hence satisfy the condition. Moreover, it leaves the disk $D_1$ on the tubular neighborhood of the second branch, on its left. See Figure \ref{figu:level2} for a representation with the four branches in a same picture.

To determine how to end the branches in $D_2 - D_1$, we first divide the Cantor set included in that annulus into $16$ chunkes. Each branch will be associated to $4$ chunks: we end them as in Figure \ref{figu:level2}.

\begin{figure}[!h]
\centering
\labellist
\small\hair 2pt
\pinlabel $D_2$ at 20 80
\endlabellist
\centering
\includegraphics[scale=0.8]{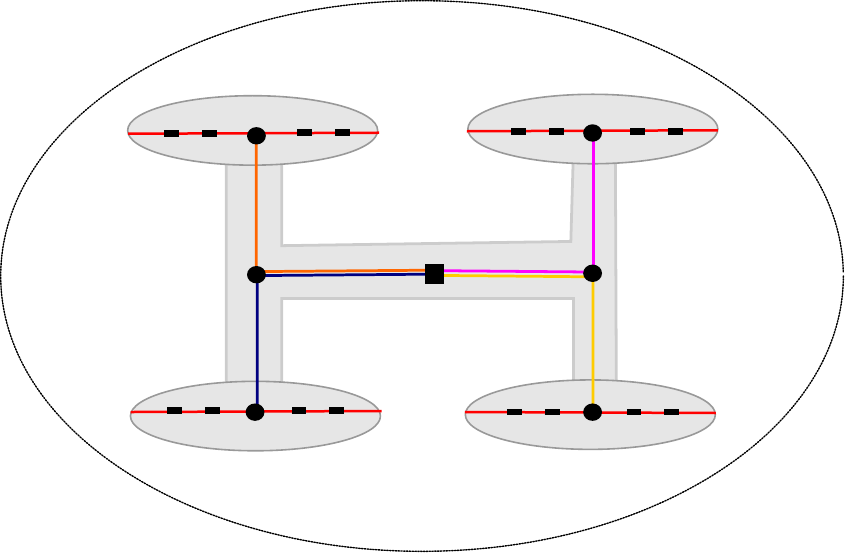}
\caption{Tubular neighborhood of level $2$ and the four disks}
\label{figu:level2-6}
\end{figure}

Note that the level $2$ of the tree is embedded in such a way so that:
\begin{itemize}
\item It is included in $D_2$.
\item All the vertices are included in the equator.
\item Each ray disjoint from level $2$ begins like level $1$.
\item We can choose a tubular neighborhood of the embedded second level of the tree and $2^2$ closed disks included in $D_2 - D_1$ so that the situation is as in Figure \ref{figu:level2-6}:
\begin{itemize}
\item Each disk contains one vertex at distance $2$ from the root and some points of the Cantor set.
\item Each boundary of the four disks intersect the equator exactly twice.
\item If we consider the intersection between one of the four disks and the equator, the vertex splits the Cantor set into two parts along this intersection. 
\end{itemize}
\end{itemize}

\section{From level $k$ to level $k+1$}

We now describe how to extend the embedding of the level $k$ for any $k$ to the embedding of the level $k+1$. This is a generalization of the construction of level $2$ from level $1$.

\begin{figure}[!h]
\centering
\labellist
\small\hair 2pt
\endlabellist
\centering
\includegraphics[scale=0.65]{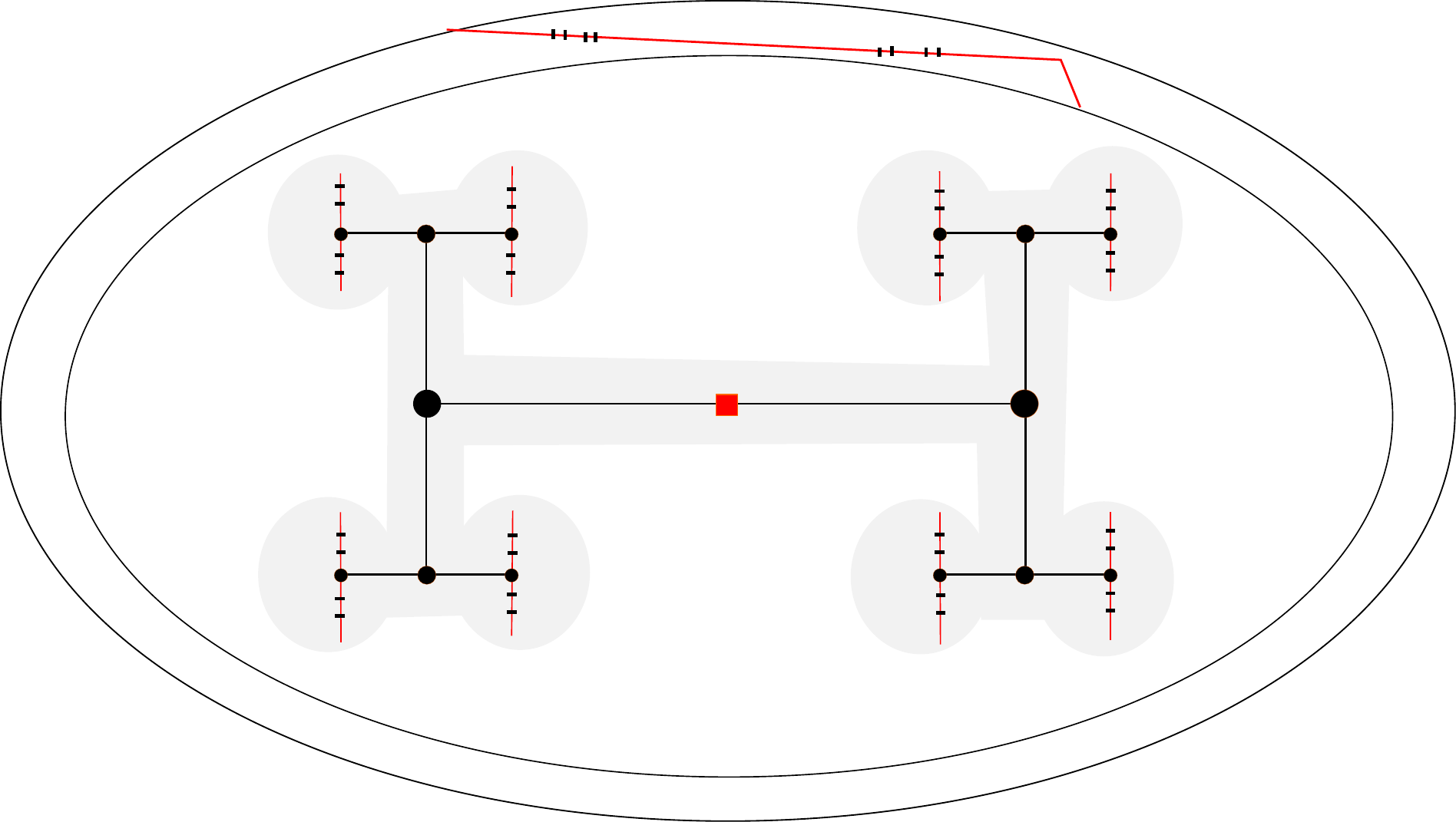}
\caption{Tubular neighborhood of level $3$ and the $2^3$ disks}
\label{figu:levelk-1}
\end{figure}

We assume that the level $k$ of the tree is embedded in such a way so that property $P(k)$ is satisfied:
\begin{itemize}
\item It is included in $D_k$.
\item All the vertices are included in the equator.
\item Each ray disjoint from level $k$ begins like level $k-1$.
\item We can choose a tubular neighborhood of the embedded $k^{th}$ level of the tree and $2^k$ closed disks included in $D_k -D_{k-1}$ so that (see Figure \ref{figu:levelk-1} for an other example with $k=3$):
\begin{itemize}
\item Each disk contains one vertex at distance $k$ from the root and some points of the Cantor set.
\item Each boundary of the $2^k$ disks intersect the equator exactly twice.
\item If we consider the intersection between one of the $2^k$ disks and the equator, the vertex splits the Cantor set into two parts along this intersection. 
\end{itemize}
\end{itemize}

\begin{figure}[!h]
\centering
\labellist
\small\hair 2pt
\endlabellist
\centering
\includegraphics[scale=0.65]{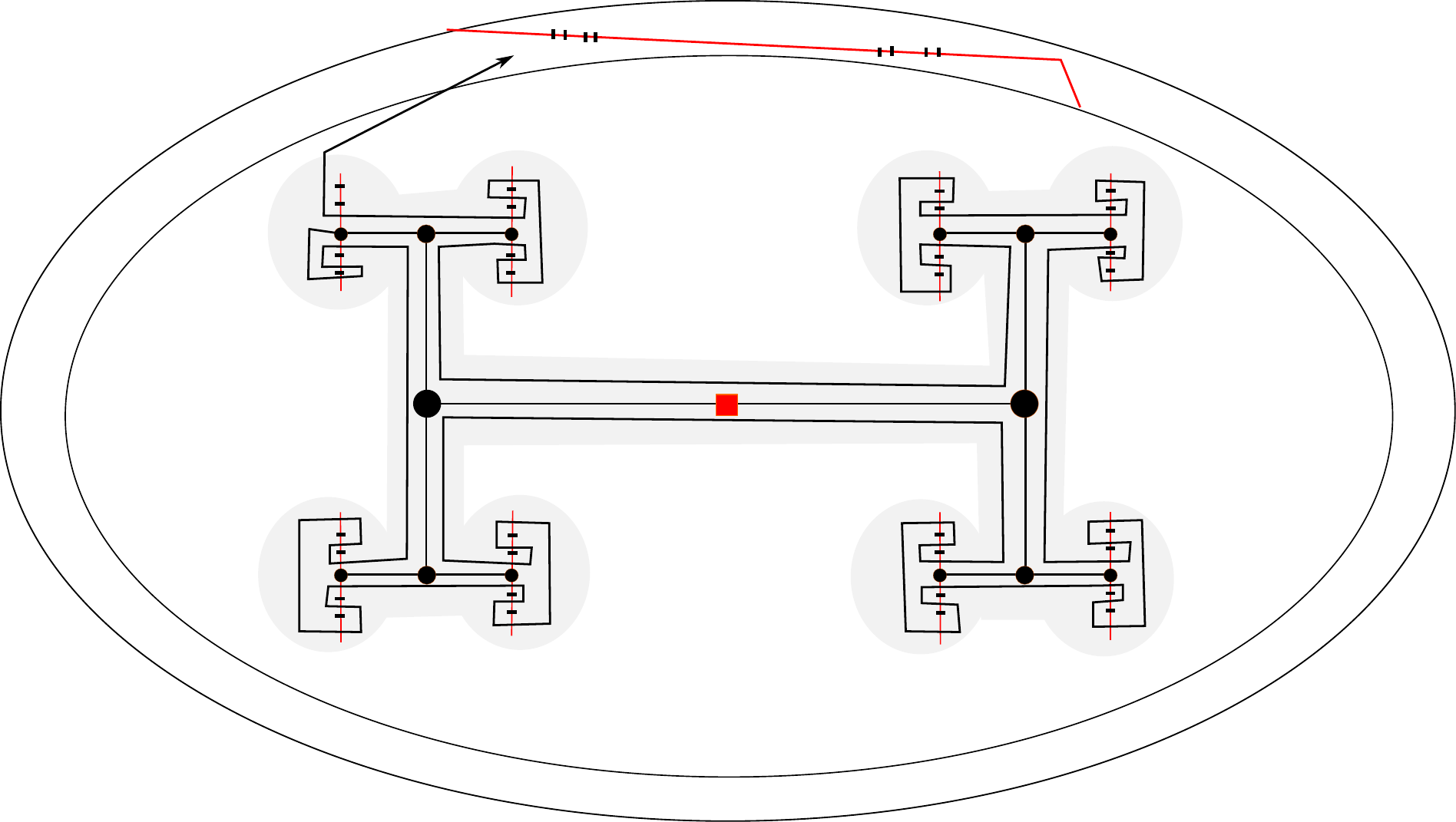}
\caption{First branch of level $k+1$}
\label{figu:levelk-2}
\end{figure}

We draw a branch of level $k+1$ in such a way so that this branch contains arcs homotopic to the boundary of the tubular neighborhood of level $k$, as in Figure \ref{figu:levelk-2}. This implies in particular that any ray disjoint from this branch has to begins like level $k$. Note that the pattern used to draw this branch in Figure \ref{figu:levelk-2} does not depend on $k$: we just need to draw a mushroom around each vertex of level $k$, which is always possible because we assumed that the Cantor set and the equator were as we needed.

We then draw the second branch, which is the one on the left of the first one in the planar representation of the tree. We draw it has we did for level $2$ in Figure \ref{figu:level2-4}: we turn around an accessible chunk of Cantor set on is right, and then go back to infinity to arrive on the left of the first branch. We then follow the second branch in a tubular neighborhood until it leaves the disk $D_k$. We draw all the branches of level $k+1$ this way, as we did for level $2$.

They all have the property that any ray disjoint from one of them has to begin like a branch of level $k$. Moreover, we can extend them in $D_{k+1} - D_k$ and attach them to the equator in such a way so that property $P(k+1)$ is satisfied.

\section{Embedded tree}
By induction, we can embed the whole tree. By construction, any ray disjoint from one of its branches has to begins like one of its branches. In other words, the embedding represents a clique of a Cantor set of disjoint high-filling rays.

\end{document}